\documentclass[12pt]{article}
\usepackage{amsmath,amsthm,amsfonts,amssymb,mathrsfs,bm,graphicx,amscd,epsfig,psfrag,verbatim,hyperref}
\usepackage[usenames]{color}
\usepackage{fancyhdr}
\parskip 	\bigskipamount

\newtheorem{remark}{Remark}[section]

\def\Z{\mathbb{Z}}

\def\R{\mathbb{R}}

\def\O{\Omega}

\def\d{\delta}
\def\del{\delta}

\def\eps{\epsilon}

\def\b{\mathcal{B}}

\def\ph{\varphi}

\def\var{\mathrm{Var}}

\def\c{\nka}

\def \d{\mathrm{d}}

\def\tr{\mathrm{tr}}
\def\eps{\varepsilon}
\def\O{\Omega}

\def\b{V_{\eps}}

\def \del{\delta}

\def \Z{\boldsymbol{Z}}

\def\L{\Lambda}

\def\rt{\rho_\tr^{(2)}}

\def\ph{\varphi}
\def\c{\complement}
\def\eq{\mathrm{eq}}
\def\eps{\varepsilon}
\def\O{\Omega}
\def \del{\delta}

\def\b{\beta}
\def\ex{\mathbf{x}}
\def\ey{\mathbf{y}}
\def\ez{\mathbf{z}}

\renewcommand{\l}[0]{\left }
\renewcommand{\r}[0]{\right}

\makeatletter
\renewcommand*{\@cite@ofmt}{\hbox}
\makeatother

\begin{document}

\title{Number rigidity in superhomogeneous random point fields}
\author{
\begin{tabular}{c}
{Subhro Ghosh}\\  
\end{tabular}
\and
\begin{tabular}{c}
{Joel Lebowitz}
\end{tabular}
}
\date{\today}
\maketitle

\begin{abstract}
We give sufficient conditions for the number rigidity of a translation invariant or periodic point process on $\R^d$, where $d=1,2$. That is, the probability distribution of the number of particles in a bounded domain $\L \subset \R^d$, conditional on the configuration on $\L^\c$, is concentrated on a single integer $N_\L$. These conditions are : (a) the variance of the  number of particles in a bounded domain $\O \subset \R^d$ grows slower than the volume of $\O$ (a.k.a. superhomogeneous  point processes), when $\O \uparrow \R^d$ (in a self-similar manner), and (b) the truncated pair correlation function is bounded by $C_1[|x-y|+1]^{-2}$ in $d=1$ and by $C_2[|x-y|+1]^{-(4+\eps)}$ in $d=2$. These conditions are satisfied by all known processes with number rigidity [\cite{GP},\cite{G},\cite{PS},\cite{AM},\cite{Bu},\cite{BuDQ}], \cite{BBNY}, and many more. in $d=1,2$. We also observe, in the light of the results of \cite{PS}, that no such criteria exist in $d>2$.
\end{abstract}

\section{Introduction}

Point processes on $\R^d$ (or $\Z^d$) are measures $\mu(\d X)$ giving rise to consistent probability measures $\mu_\L$ on
configurations of particles $X_{\L}$ on all regions $\L \subset \R^d$ (or $\Z^d$). When the volume of $\L$, denoted $|\L|$, is finite, this probability is concentrated on configurations $X_\L=\{x_1,\cdots,x_N; x_l \in \L\}$, $N$ finite. From $\mu (\d X)$, we can find the conditional probabilities $\mu(\d X_\L | \L^\c)$ of configurations in $\L$, given the configuration in $\L^\c = \R^d \setminus \L$. We shall always assume that $\mu$ is ergodic and either translation invariant or periodic.

%

The simplest example of a point process on $\R^d$ is the translation invariant Poisson process with density (intensity) $\rho>0$. For this process, there is no information gained about $X_\L$ from knowing about $X_{\L^\c}$, so $\mu(\d X_\L | X_{\L^\c})=\mu(\d X_\L)$. In general, for equilibrium systems with sufficiently rapidly decaying (tempered) interaction potentials $U(X)$, the infinite volume Gibbs measures, $\mu_\eq$, describe a point process whose conditional probabilities satisfy the DLR (Dobrushin, Lanford, Ruelle) equation 
\begin{equation} 
\label{DLR}
 \mu_\eq (\d X_\L | X_{\L^\c})  = \exp [ - \b U(X_\L | X_{\L^\c})  ] \bigg/ \l( \int_\L \exp(-\b U(X_\L | X_{\L^\c})) \d X_\L\r),
\end{equation}
where $\b$ is the inverse temperature \cite{R},\cite{Ge}. This gives a probability distribution for the number of particles in $\L$ given $X_{\L^\c}$, which increases as $|\L|$ increases. 

In contrast to the above situation, we shall be interested here in the case, called ``number rigidity'' by \cite{GP}, where $X_{\L^\c}$ determines the precise number of particles in $\L$. An early example of such a property (presumably the earliest) was proven by Aizenman and Martin (\cite{AM}) for one dimensional Coulomb systems in which the interactions are very long range. In this case, $U$ is a sum of interactions between particles which increases linearly with distance. For such systems the DLR equations are not well defined for the limiting infinite volume Gibbs measure. \cite{AM} considered two cases of charge neutral systems: \begin{enumerate}
\item A one component plasma (a.k.a. jellium) in which there is a uniform background of negative charge with density $\rho$ and  charge one positive point particles of average density $\rho$. In this case the extremal states are periodic with period $\rho^{-1}$.
\item A multicomponent charged system with both positively and negatively charged particles, for which $\mu(\d X)$ is translation invariant.
\end{enumerate}
In the second case, number rigidity is replaced by net charge rigidity, i.e. $X_{\L^\c}$ determines the total net charge in $\L$. For the OCP, net charge is just $N_\L - \langle N_\L \rangle = N_\L - \rho |\L|$.

The system considered by \cite{GP} is the determinantal point process corresponding to the OCP in $d=2$, at a value of the reciprocal temperature $\beta=2$. This measure is known to be identical to that of the eigenvalue distribution of the Ginibre ensemble (the weak limit of eigenvalues of matrices with i.i.d. complex Gaussian entries with mean 0 and variance 1), \cite{Gi-2}. For the  2D OCP (or the multicomponent Coulomb system) the interaction  between the charges grows like the logarithm of the distance between them. The model considered by \cite{G} is also a determinantal point process corresponding to the eigenvalue distribution of the Gaussian Unitary Ensemble (GUE) or the Circular Unitary Ensemble (CUE).  This distribution is the same as that of the equilibrim 1-D ``Dyson logarithmic gas'' (DLG), again at $\b=2$ (see \cite{For}). This system also has a uniform negative charge background, but the interaction between the positive charges only grows logarithmically, rather than linearly, with distance. 
Other examples of ``number rigidity'' were described by \cite{PS}, \cite{Bu}, \cite{BuDQ} and \cite{BBNY}.

In this note we give sufficient conditions for rigidity which covers the above cases, as well as the \cite{PS} case in $d=1,2$, and extends to a large class of other point processes in $d=1,2$, e.g. the G-process and related self-correcting queuing processes described in \cite{GLS} as well as to the Conway-Radin tilings (\cite{Ra}).

 \section{Superhomogeneous systems}

 In this section, we show that for translation invariant point processes having sub-volume growth of number variance (i.e., variance of the number of particles in a ball of radius $R$ is $o(R^d)$), the variance of linear statistics assumes a particularly simple form.  
 
 To this end, we consider the truncated total pair correlation function (\cite{MY}) \[G(x,y)=\rho_1(x) \del(x,y) + \rt(x,y),\] where  $\rt=\rho_2(x,y)-\rho_1(x)\rho_1(y)$ is the truncated pair correlation function, $\rho_1(x)$ is the one particle density, and $\del(x,y)$ is the Dirac (Kronecker) delta function. Clearly, we have $G(x,y)=G(y,x)$. Moreover, for a translation invariant process, $G(x,y)=G(x-y)$. Consider now a sequence of domains $\L \uparrow \R^d$ in a regular (e.g. self-similar) way. In this setting, we have (see \cite{MY})
 \[\var(N_\L)\equiv \int_\L \int_\L G(x,y) \d x \d y= |\L|\int_{\R^d} G(x,y)\d y + o( |\L|). \] For processes with sub-volume growth of the variance of $N_\L$, this  gives \begin{equation} \label{symm} \int_{\R^d} G(x,y)\d x = \int_{\R^d} G(x,y)\d y =0. \end{equation} 

Note that for point processes satisying the DLR equation (and some other mild conditions), the variance of the number of particles in $\L$ (denoted $N_\L$) is bounded below by $c|\L|, c>0$, for $\b$ finite (\cite{Gi-1},\cite{LPRS}). Such processes are therefore not superhomogeneous and also not rigid for finite $\b$.

 \textbf{Linear Statistics}
 
 For a test function $\ph$ and a point process $\Pi$, we define the linear statistic \[\int \ph \d [\Pi]= \sum_{x\in \Pi} \ph(x).\]
 We then have 
 \begin{equation} \label{vls1} \var\l( \int \ph \d [\Pi]  \r) = \int \int \ph(x) \ph(y) G(x,y)\d x \d y. \end{equation}
 
 But for systems satisfying \eqref{symm} we have  $\int |\ph(x)|^2 G(x,y) \d x \d y =0$. This implies that \eqref{vls1} can be rewritten as   \begin{equation} \label{vls2} \var\l( \int \ph \d [\Pi]  \r) = - \frac{1}{2}\int \int |\ph(x) - \ph(y)|^2 G(x,y)\d x \d y. \end{equation} However, due the presence of $ |\ph(x) - \ph(y)|^2$, the $\del(x-y)$ component of $G(x,y)$ contribute 0 to the integral on the right hand side of \eqref{vls2}. This enables us to write  \begin{equation} \label{vls3} \var\l( \int \ph \d [\Pi]  \r) = - \frac{1}{2}\int \int |\ph(x) - \ph(y)|^2 \rt(x,y)\d x \d y. \end{equation}

 Our goal now is to construct,for any $\eps >0$, a function $\ph^\eps$ such that $\ph^\eps \equiv 1$ on the unit ball in $\R^d,d=1,2$, and $\var \l( \int \ph^\eps \d [\Pi] \r) <\eps$. Using the general strategy in \cite{GP} or \cite{G}, we can then deduce number rigidity of $\Pi$.

 \section{Variance of linear statistics in $d=2$} 
  
 To illustrate our method, we consider first the case where \[|\rt(x,y)|\le C \exp(-\gamma|x-y|).\]
 
 We begin with a non-negative $C_c^2$ function $\Phi$ supported on a disk of radius $K$ such that $\|\Phi\|_\infty \le 1$ and $\Phi \equiv 1$ on the unit disk.
 
 We set $\Phi_R(x)=\Phi(x/R)$. The variance of the linear statistic corresponding to $\Phi_R$ has the bound 
 \[ \var\l( \int \Phi_R \d [\Pi]\r) \le C \int \int \exp(-\gamma|x-y|) |\Phi_R(x) - \Phi_R(y)|^2 dm(x)dm(y) \] where $\d m(x)$ is the Lebesgue measure on $\R^2$.
 We upper bound the above by integrals over two regions, the point being that due to the support properties of $\Phi$ the integrand vanishes outside $A_1 \cup A_2$ defined below:
\[A_1:=\{|x| \le KR  \},\]
\[A_2:=\{|y| \le KR\}.\]
 Due to symmetry between $x$ and $y$, it suffices to estimate from above
 \[ \int \int_{A_1} \exp(-\gamma|x-y|) |\Phi_R(x) - \Phi_R(y)|^2 dm(x)dm(y).\] 
%
 For the integral over $A_1$, we expand $\Phi_R$ in a Taylor series in $y$ as \[ \Phi_R(y)=\Phi_R(x)+ \frac{ \langle \nabla \Phi(\frac{x}{R}) , y-x \rangle }{R} + h(x,y) , \] where the error term $h(x,y)$ is bounded as \[ |h(x,y)| \le \frac{C_2(\Phi)|y-x|^2}{R^2}, \] and the quantity $C_2(\Phi)$ is given by $C_2(\Phi)=A \sup_z \|D^2\Phi(z)\|_{2}$ for some universal constant $A$; in other words the supremum of the 2-norm of the Hessian matrix $D^2\Phi(z)$ of $\Phi$. Analogously, define $C_1(\Phi)=A \sup_z \|\nabla \Phi (z)\|_{2}$. The upshot of this is that \[|\Phi_R(x) - \Phi_R(y)|^2 \le \|\nabla \Phi (\frac{x}{R})\|^2 \frac{|y-x|^2}{R^2} + \frac{2C_1(\Phi)C_2(\Phi)|y-x|^3}{R^3} + \frac{C_2(\Phi)^2|y-x|^4}{R^4}.  \]
 
 Observe that 
 \begin{align*} &\int \int_{|x| \le KR} \frac{2C_1(\Phi)C_2(\Phi)|y-x|^3}{R^3}\exp(-\gamma|y-x|)dm(y)dm(x) \\  \le  &\frac{2C_1(\Phi)C_2(\Phi)}{R^3}
  \l(\int_{|x|\le KR}dm(x)\r)\l( \int |y-x|^3\exp(-\gamma|y-x|)dm(y) \r) \\ &= K^2C_1(\Phi)C_2(\Phi)C_3(\gamma)/R. \end{align*}
 Similarly, we can deduce that 
 \[ \int_{|x|\le KR} \frac{C_2(\Phi)^2|y-x|^4}{R^4}\exp(-\gamma|x-y|)dm(y)dm(x) \le K^2C_2(\Phi)^2C_4(\gamma)/R^2. \]
 We are thus left with the term \begin{align*} &\int \int_{|x|\le KR}  \|\nabla \Phi (\frac{x}{R})\|^2 |y-x|^2 \exp(-\gamma|y-x|) dm(y)dm(x) \\ = &C_5(\gamma)\int_{|x|\le KR}\frac{|\nabla \Phi(x/R)|^2}{R^2}dm(x) \\ = &C_5(\gamma) \int_{|u|\le K} |\nabla \Phi(u)|^2 dm(u) \\  = &C_5(\gamma) \|\nabla \Phi\|_2^2. \end{align*}
 
 
 Thus, for $\Phi$ as described above, we have the estimate 
 \[\var\l( \int \Phi_R \d [\Pi]\r) \le 2C_5(\gamma)\|\nabla \Phi\|_2^2 + 2K^2C_1(\Phi)C_2(\Phi)C_3(\gamma)/R  +  2K^2C_2(\Phi)^2C_4(\gamma)/R^2.\] where $C$ is a universal constant.
 
 Now we select $\Phi$ such that $\|\nabla \Phi\|_2^2 \le \eps/6C_5(\gamma)$ (this can be done as in \cite{GP} for the case of the Ginibre ensemble).  Depending on $\Phi$, we choose $R$ so large that $\max\{2K^2C_1(\Phi)C_2(\Phi)C_3(\gamma)/R ,  2K^2C_2(\Phi)^2C_4(\gamma)/R^2\} \le \eps/3$. For such a choice of $\Phi$ and $R$, we can take $\phi^\eps=\Phi_R$ and we shall have $\var\l( \int \phi^\eps \d [\Pi]\r) <\eps$, as desired.
 
%

\textbf{Power-law decay of correlations}

 In this section, we present a refinement of the argument in the previous section, which implies a similar variance bound (and hence, rigidity of numbers) when we have a power law decay of correlations. We will see that a power law decay of the truncated pair correlation function \begin{equation} \label{powerlaw} |\rt(x,y)|\le C(1+|x-y|^{4+\eps})^{-1}\end{equation} will suffice. 
 
 To this end, we once again consider the variance of linear statistics, and consider a non-negative $C_c^2$ function $\Phi$ supported on a disk of radius $K$ such that $\|\Phi\|_\infty \le 1$ and $\Phi \equiv 1$ on the unit disk. With notations as before, the variance is given by  \[ \var\l( \int \Phi_R \d [\Pi]\r) \le C \int \int \rt(x,y) |\Phi_R(x) - \Phi_R(y)|^2 dm(x)dm(y). \] As before, we consider the integral only in the region $\{(x,y) \in A_1 \cup A_2 \}$, where $A_1:=\{|x| \le KR  \}, A_2:=\{|y| \le KR\}$.
 
 By symmetry, it suffices to obtain an upper bound on the integral only over the region  $A_1$. We will work with the decomposition $A_1=(A_1\cap A_2) \cup (A_1 \cap A_2^\c) $.

 For integration over $A_1\cap A_2$, as in the previous section, we start with the estimate   \begin{equation}\label{est1}|\Phi_R(x) - \Phi_R(y)|^2 \le \|\nabla \Phi (\frac{x}{R})\|^2 \frac{|y-x|^2}{R^2} + \frac{2C_1(\Phi)C_2(\Phi)|y-x|^3}{R^3} + \frac{C_2(\Phi)^2|y-x|^4}{R^4}.  \end{equation}
 Of the various quantities on the right hand side of \eqref{est1}, we estimate \[\int \int_{A_1 \cap A_2} \Phi (\frac{x}{R})\|^2 \frac{|y-x|^2}{R^2} \rt(x,y) \d m(x) \d m(y)  \] from above by \[ \int_{|x| \le KR} \l( \int |y-x|^2\rt(x,y) \d m(y)\r) \|\nabla \Phi (\frac{x}{R})\|^2 \frac{\d m(x)}{R^2} = B  \|\nabla \Phi \|_2^2, \] where $B=\int |y-x|^2\rt(x,y) \d m(y) <\infty$. It is easy to see that $B$ is independent of  $x$, and its finiteness follows from the power law decay \eqref{powerlaw} assumed on $\rt$.
 
 Our task now boils down to showing that $\frac{1}{R^3} \int_{A_1 \cap A_2 } {|y-x|^3} |\rt(x,y)|\d m(x) \d m(y)$ and $\frac{1}{R^4} \int_{A_1 \cap A_2} {|y-x|^4} |\rt(x,y)| \d m(x) \d m(y) $ are $o(1)$ as $R \to \infty$. We will show this in detail for the first integral; the argument is similar for the second one. 
 
 To upper bound $\frac{1}{R^3} \int_{A_1 \cap A_2} {|y-x|^3} |\rt(x,y)|\d m(x) \d m(y)$, we proceed as:
 \[ \frac{1}{R^3} \int_{A_1 \cap A_2} {|y-x|^3} |\rt(x,y)| \d m(y) \d m(x) \le  \frac{C}{R^3} \int_{|x|\le K R} R^{1-\eps} \d m(x) = o(1),   \] where in the last inequality we have used \eqref{powerlaw}.
 
 This leaves us with the integral over $A_1 \cap A_2^\c$.  To handle this, we notice that on  $A_1 \cap A_2^\c$, we have $ \Phi_R(y)=0$ due to the support properties of $\Phi$. Consequently, \begin{align*} &\int\int_{A_1 \cap A_2^\c}|\Phi_R(x) - \Phi_R(y)|^2 |\rt(x,y)| \d m(x) \d m(y) \\ = &\int_{A_1} |\Phi_R(x)|^2 \l( |\rt(x,y)| \d m(y)  \r) \d m(x) \\ \le &  C \int_{|x| \le K R} |\Phi_R(x)|^2 \cdot \frac{1}{R^{2+\eps}} \d m(x) \\ =  &\frac{1}{R^\eps}\int |\Phi_R(x)|^2 \frac{\d m(x)}{R^2} \\ = &\frac{1}{R^\eps} \|\Phi\|_2^2 \\ = &o(1). \end{align*}

  \begin{remark}
  A careful accounting in the above argument shows that, in fact, \[\lim_{R\to \infty}\var\l( \int \Phi_R \d [\Pi]\r)=C\|\nabla \Phi\|_2^2\] for point processes with $4+\eps$ (or faster) correlation decay.
 \end{remark}
  
  \begin{remark}
   It is also clear that a power law $4+\eps$ is not strictly necessary, a milder growth faster than $R^{-4}$, like $R^{-4}(\log R)^{-1}$, should also be enough for the method to work.
  \end{remark}

  \begin{remark}
   The analysis of the case $d=1$ is simpler than the case $d=2$. The method and the result is the same as in \cite{G}. All that is required for the rigidity of superhomogeneous systems in $d=1$ is $|\rt(x,y)| \le C[1+|x-y|]^{-2}$.
  \end{remark}

\section{Rigidity in $d \ge 3$} 
 

We have shown that any superhomogeneous process in $\R^d, d=1,2$, whose $|\rt|$ satisfies the power law decay above has number rigidity. 
It is clearly of interest to consider what happens in higher dimensions. In this context, it is of interest to look at the point processes given by i.i.d. perturbations of a lattice, i.e. each lattice point $z \in \R^d$ is shifted to $z+x \in \R^d$ with a probability distribution $h(\d x)$. It has been shown by Peres and Sly (\cite{PS}) that for Guassian lattice perturbations (i.e., $h(x)=(2 \pi)^{-d}\exp(-\frac{-x^2}{2 \sigma^2})$), there is a phase transition in the rigidity behaviour in dimension $d \ge 3$. Namely, there is a critical  $\sigma_c>0$ such that if the standard deviation $\sigma$ of the perturbations satisfies $\sigma < \sigma_c$, then there is rigidity of the number of particles, whereas there is no rigidity when $\sigma > \sigma_c$. 

For i.i.d. Gaussian perturbations of $\Z^d$, it is always the case that the point process is superhomogeneous and has fast decay of correlations. This can be seen from the fact that when the density of the perturbing random variable is $h \in L^{\infty}$, the particle density $\rho_1(\ex)$ and the truncated pair correlation function $\rt(\ex,\ey)$ are periodic and are given by
\[\rho_1(\ex)=\sum_{\mathbf{z} \in \Z^d} h(\ex-\mathbf{z})\] and 
\[\rho_\tr^{(2)}(\ex,\ey)=-\sum_{\ez \in \Z^d } h(\ex-\ez)h(\ez -\ey). \] 
It follows immediately from the definition of $G$ that $\int_{\R^d} G(\ex,\ey)\d \ey =0$, implying that these systems are superhomogeneous.

When $h$ decays fast enough, $|\rt(\ex,\ey)|$ is bounded above in absolute value by $C_h h(x-y)$, where $C_h$ is a quantity that depends on the function $h$. This is in particular true in the case of Gaussian perturbations, which implies that the truncated pair correlation function decays like a Gaussian in that case. More generally this would be true when $h(\ex)$ decays faster than $C|\ex|^{-\gamma}$ with $\gamma > d$, as $|\ex| \to \infty$. The main idea is that when $\int_{\R^d} h(\ex) \d \ex $ is finite, the sum  $\sum_{\ez \in \Z^d } h(\ex-\ez)h(\ez -\ey)$ is dominated by the terms $h(\ex-\ez)h(\ez -\ey)$ when $\ez=\ex$ and $\ez=\ey$.
 
For $d \ge 3$, this provides us a concrete counterexample to any conjectured relationship between rigidity and decay of correlations, even under the assumption of superhomogeneity. However, this example is no longer valid in $d =1,2$, where it has been shown in \cite{PS} that the existence of the first and second moments of $h$ is sufficient for the point process to exhibit rigidity of numbers. Therefore, to understand the precise relationship between rigidity, decay of correlations and superhomogeneity in higher dimensions remains a delicate open question. 

\section{Concluding remarks}

All the point processes for $d=1,2$ for which rigidity has been proven satisfy our conditions. For 2-D systems like the Ginibre ensemble, the zeroes of the standard planar Gaussian analytic function on the complex plane and the processes considered by \cite{PS}, the variance of the number of points grows like the surface area of the domain. This is the slowest possible rate of growth of variance for isotropic point processes (\cite{Be}). This is in fact expected (and in many cases proven) for the variance of charge (particle number in OCP) for Coulomb systems in $d>1$, while in $d=1$, the variance is bounded, $|\partial \L|=2$ (\cite{MY},\cite{L}). For the Dyson log gas, the variance grows like $\log |\L|$. It satisfies, for $\b \le 2$, the bound $|\rt(x,y)|\le C [1+|x-y|]^{-2}$ (\cite{For}). Note that the integrability of $|y|\rt(y)$ is not required. Another family of superhomogeneous processes is given by determinantal point processes on Euclidean spaces whose kernels are projections when thought of as an integral operator. This includes the systems studied by \cite{G} and by \cite{Bu} and \cite{BuDQ}. A third family of examples of such processes is given by i.i.d. perturbations of the lattice $\Z^d$ studied by \cite{PS}. This class of examples is perhaps the one for which superhomogeneity is the easiest to verify.

We expect our results to extend to point processes on lattice systems, in $d=1,2$, under the same conditions of superhomogeneity and decay of correlations. It is an interesting question whether superhomogeneity (possibly coupled with good decay) is necessary to guarantee rigidity behaviour in a point process. A heuristic reason for such a conjecture is that when the variance grows like $|\L|$, it behaves in an additive way for two adjacent domains which seems to suggest that surface effects become negligible for large $|\L|$, which is inconsistent with number rigidity.


\section{Acknowledgements}
We thank Michael Aizenman and Peter Forrester for many helpful discussions on the rigidity of Coulomb systems. The work of J.L.L. was supported by the NSF grant DMR1104500.

\begin{tabular}{l r}
\textsc{Subhroshekhar Ghosh}  \hspace{100 pt} & \textsc{Joel L. Lebowitz} \\
\textsc{Dept of ORFE}  & \textsc{Depts of Mathematics and Physics} \\
\textsc{Princeton University} & \textsc{Rutgers University} \\
email: subhrowork@gmail.com & email: lebowitz@math.rutgers.edu \\ 
\end{tabular}

 \end{document}